\documentclass[12pt]{article}
\usepackage{amsmath,epsfig,subfigure,amsthm,amsfonts,amsbsy,amssymb,latexsym,amsxtra}
\usepackage[usenames]{color}
\usepackage{hyperref}

\usepackage{pst-infixplot}
\usepackage{pst-plot}
\usepackage{pstricks}

\usepackage[sort]{natbib}
\makeatletter \@addtoreset{equation}{section}

\newtheorem{thm}{Theorem}[section]
\newtheorem{prop}[thm]{Proposition}
\newtheorem{lem}[thm]{Lemma}

\theoremstyle{definition}

\newtheorem{defn}[thm]{Definition}

\newcommand{\thmref}[1]{Theorem~{\rm \ref{#1}}}
\newcommand{\lemref}[1]{Lemma~{\rm \ref{#1}}}

\newcommand{\propref}[1]{Proposition~{\rm \ref{#1}}}

\newcommand{\beq}{\begin{equation}}
\newcommand{\eeq}{\end{equation}}
\newcommand{\bed}{\begin{displaymath}}
\newcommand{\eed}{\end{displaymath}}
\newcommand{\ben}{\begin{eqnarray*}}
\newcommand{\een}{\end{eqnarray*}}

\newcommand{\bea}{\bed\begin{array}{rl}}
\newcommand{\eea}{\end{array}\eed}
\newcommand{\ad}{&\!\!\!\disp}
\newcommand{\aad}{&\disp}
\newcommand{\barray}{\begin{array}{ll}}
\newcommand{\earray}{\end{array}}

\def\F{{\cal F}}
\def\op{{A}}

\def\rr{{\mathbb R}}

\newcommand{\e}{\varepsilon}

\newcommand{\R}{\mathcal R}
\newcommand{\ex}{{\mathbf E}}
\newcommand{\taub}{{\tau_{b^*}}}

\newcommand{\xz}{\ensuremath{x_0}}

\def\({\left(}
\def\){\right)}

\newcommand{\nd}{\noindent}

\newcommand{\disp}{\displaystyle}

\def\one{{\hbox{1{\kern -0.35em}1}}}

\newcommand{\set}[1]{\left\{#1\right\}}

\title{On Optimal Harvesting in Stochastic Environments: Optimal Policies in a Relaxed Model}
\author{
     Richard H. Stockbridge\thanks{Department of Mathematical Sciences, University of
    Wisconsin-Milwaukee, Milwaukee, WI 53201, {\tt  stockbri@uwm.edu}.
    The research of this author was supported in part by the U.S. National Security Agency under Grant Agreement Number H98230-09-1-0002.  The United States Government is authorized to reproduce and distribute reprints
notwithstanding any copyright notation herein.}
     \and Chao   Zhu\thanks{Department of Mathematical Sciences, University of
    Wisconsin-Milwaukee, Milwaukee, WI 53201, {\tt zhu@uwm.edu}. The research of this author was supported in part by a grant from the  UWM Research
    Growth Initiative.} }

\begin{document}

\maketitle

\begin{abstract}
This paper examines the objective of optimally harvesting a single species in a stochastic environment.  This problem has previously been analyzed in \cite{Alvarez} using dynamic programming techniques and, due to the natural payoff structure of the price rate function (the price decreases as the population increases), no optimal harvesting policy exists.  This paper establishes a relaxed formulation of the harvesting model in such a manner that existence of an optimal relaxed harvesting policy can not only be proven but also identified.  The analysis imbeds the harvesting problem in an infinite-dimensional linear program over a space of occupation measures in which the initial position enters as a parameter and then analyzes an auxiliary problem having fewer constraints.  In this manner upper bounds are determined for the optimal value (with the given initial position); these bounds depend on the relation of the initial population size to a specific target size.  The more interesting case occurs when the initial population exceeds this target size; a new argument is required to obtain a sharp upper bound.  Though the initial population size only enters as a parameter, the value is determined in a closed-form functional expression of this parameter.

\medskip
\nd {\bf Key Words.}
Singular stochastic control,
linear programming, relaxed control.

\medskip
\nd {\bf AMS subject classification.} 93E20, 60J60.

\end{abstract}

\setlength{\baselineskip}{0.23in}
\section{Introduction}\label{sect-intro}

This paper examines the problem of optimally harvesting a single species that lives in a random environment.  Let $X$ be the process denoting the size of the population and
 $Z$ denote the cumulative amount of the species harvested.  We assume $X(0-)=x_0 >0$, $Z(0-)=0$, and $X$ and $Z$ satisfy
\begin{equation}\label{harvested-pop}
dX (t) = b(X(t) )dt + \sigma(X(t) )dW(t) - dZ(t),
\end{equation}
in which $W(\cdot)$ is a $1$-dimensional standard Brownian motion that provides the random fluctuations in the population's size, and $b$ and $\sigma$ are real-valued functions.  We assume that $b$ and $\sigma$ are such that in the absence of harvesting the population process $X$ takes values in $\rr_+$ and that $\infty$ is a natural boundary so that the population will not explode to $\infty$ in finite time.  The boundary $0$ may be an exit or a natural boundary point but may not be an entrance point; this indicates that the species will not reappear following extinction.  Note that $X(0)$ may not equal $X(0-)$ due to an instantaneous harvest $Z(0)$ at time $0$ and the process $Z$ is restricted so that $\Delta Z(t) := Z(t) - Z(t-) \leq X(t-)$ for all $t \geq 0$.  This latter condition indicates that one cannot harvest more of the species than exists.  Let $r>0$ denote the discount rate and $f$ denote the marginal yield for harvesting.  The objective is to select a harvesting strategy $Z$ so as to maximize the expected discounted revenue
\begin{equation} \label{income-J-defn}
J(x_0, Z) :=\ex_{\xz} \left[\int_{0}^{\tau} e^{-r s} f(X(s-) ) dZ(s)\right],
\end{equation}
where $\tau=\inf\set{t\ge 0: X(t) = 0}$ denotes the extinction time of the species.

 As a result of developments in stochastic analysis and stochastic control techniques, there has been a resurgent interest in determining the optimal harvesting strategies in the presence of stochastic fluctuations (see, e.g., \cite{A-Shepp, Brauman-02,Jorgensen-Y-96,Lungu-O,L-Oksendal,Ryan-Hanson}).  In particular, \cite{Alvarez} examines the current problem using dynamic programming techniques and determines the value function.  The paper indicates the lack of an optimal policy in the admissible class of (strict) harvesting policies by commenting that a ``chattering'' policy will be optimal.  The problem of optimal harvesting of a single species in a random environment is also studied in \cite{Song-S-Z} in which the model is extended to regime-switching diffusions so as to capture different dynamics such as for drought and non-drought conditions.  The paper also adopts a dynamic programming solution approach to determine the value function while at the same time exhibiting $\epsilon$-optimal harvesting policies since, as in the static environment of \cite{Alvarez}, no optimal harvesting policy exists.  In light of the complexities of the regime-switching model, it further identifies a condition under which the value function is shown to be continuous and a viscosity solution to the variational inequality.  

The focus of this paper is on developing a relaxed formulation for the harvesting problem under which an optimal harvesting control exists and on establishing optimality using a linear programming formulation instead of dynamic programming.  In addition, it is sufficient to have a weak solution to (\ref{harvested-pop}) rather than placing Lipschitz and polynomial growth conditions on the coefficients $b$ and $\sigma$ that guarantee existence of a strong solution.  Intuitively, relaxation completes the space of admissible harvesting rules by allowing measure-valued policies.  A benefit of the linear programming solution methodology is the analysis concentrates on the optimal value for a single, fixed initial condition, rather than seeking the value {\em function} and thus no smoothness properties need to be established about the value as a function of the initial position.

To set the stage for the relaxed singular control formulation of the model, let ${\cal D} = C^2_c(\rr_+)$ and for a function $g \in {\cal D}$, define the operators $A$ and $B$ by
\begin{align} \label{genA}
&Ag(x ) = \frac{1}{2} \sigma^2(x ) g''(x ) + b(x ) g'(x ), \mbox{ and } \\ \label{genB}
& Bg(x, z) = \begin{cases}\frac{g(x-z )-g(x )}{z}, & \text{ if } z> 0, \\
  -g'(x), & \text{ if } z=0,\end{cases}
\end{align}
where $x,z \in \rr_+$.  It\^{o}'s formula then implies
\begin{eqnarray*}
g(X(t) ) &=& g(x_0 ) + \int_0^t Ag(X(s) )\, ds + \int_0^t Bg(X(s), \Delta Z(s))\, dZ(s) \\
& & + \int_0^t \sigma(X(s) ) g'(X(s))\, dW(s), \ \ \forall g\in \mathcal D.
\end{eqnarray*}
It therefore follows that for any $g\in \mathcal D$
\begin{equation}  \label{mg}
g(X(t) ) - g(x_0 ) - \int_0^t Ag(X(s) )\, ds - \int_0^t Bg(X(s) ,\Delta Z(s))\, dZ(s)
\end{equation}
is a mean $0$ martingale.  In fact, requiring (\ref{mg}) to be a martingale for a sufficiently large collection of functions $g$ is a way to characterize the processes $(X, Z)$ which satisfy (\ref{harvested-pop}).  We turn now to a precise formulation of the model in which the processes are relaxed solutions of a controlled martingale problem for the operators $(A,B)$.

\subsection{Formulation of the Relaxed Model} \label{form}
For a complete and  separable  metric space $S$, we define $M(S)$ to be the space of Borel measurable functions on $S$, $B(S)$ to be the space of bounded, measurable functions on $S$, $C(S)$ to be the space of continuous functions on $S$, $\overline C(S)$ to be the space of bounded, continuous functions on $S$, ${\cal M}(S)$ to be the space of finite Borel measures on $S$, and ${\cal P}(S)$ to be the space of probability measures on $S$.  ${\cal M}(S)$ and ${\cal P}(S)$ are topologized by weak convergence.

Recall, the amount of harvesting is limited by the size of the population.  Define ${\cal R} = \{(x,z): 0 \leq z \leq x, x \geq 0\}$; ${\cal R}$ denotes the space on which the paired process $(X,Z)$ evolves when considering solutions of (\ref{harvested-pop}).

The formulation of the population model in the presence of ``relaxed'' harvesting policies adapts the relaxed formulation for singular controls given in \cite{kurt:01} to the particulars of the harvesting problem.  This adaptation sets the state space $E$ to be $\rr_+ $ and the control space $U = \rr_+$, with ${\cal U} = {\cal R} \subset \rr_+ \times \rr_+$.  We begin by specifying the space of measures for the relaxed harvesting policies.  Let ${\cal L}_t({\cal R})={\cal M}({\cal R} \times [0,t])$.  Define ${\cal L}({\cal R})$ to be the space of measures $\xi$ on ${\cal R} \times [0,\infty )$ such that $\xi ({\cal R} \times [0,t])<\infty$, for each $t$, and topologized so that $\xi_n\rightarrow\xi$ if and only if $\int f\, d\xi_n \rightarrow \int fd\, \xi$, for every $f\in\overline C({\cal R} \times [0,\infty ))$ with supp$(f)\subset {\cal R} \times [0,t_f]$ for some $t_f<\infty$.  Let $\xi_t\in {\cal L}_t({\cal R} )$ denote the restriction of $\xi$ to ${\cal R} \times [0,t]$.  Note that a sequence $\{\xi^n\}\subset {\cal L}({\cal R} )$ converges to a $\xi\in {\cal L}({\cal R} )$ if and only if there exists a sequence $\{t_k\}$, with $t_k\rightarrow \infty$, such that, for each $t_k$, $\xi^n_{t_k}$ converges weakly to $\xi_{t_k}$, which in turn implies $\xi^n_t$ converges weakly to $\xi_t$ for each $t$ satisfying $\xi ({\cal R} \times \{t\})=0$.

Let $X $ be an $\rr_+ $-valued process and $\Gamma$ be an ${\cal L}({\cal R})$-valued random variable.  Let $\Gamma_t$ denote the restriction of $\Gamma$ to ${\cal R} \times [0,t]$.    Then $(X,\Gamma)$ is a {\em relaxed solution} of the harvesting model if there exists a filtration $\{{\cal F}_t\}$ such that $(X,\Gamma_t)$ is $\{{\cal F}_t\}$-progressively measurable, $X(0)=x_0$,   and for every $g\in {\cal D}$,
\begin{equation}\label{mgp}
g(X(t))-\int_0^t Ag(X(s))\, ds - \int_{{\cal R} \times [0,t]} Bg(x,z)\,\Gamma (dx\times dz \times ds)
\end{equation}
is an $\{{\cal F}_t\}$-martingale, in which the operators $A$ and $B$ are given by (\ref{genA}) and (\ref{genB}), respectively.
Throughout the paper we assume that a relaxed solution $(X,\Gamma)$ exists and that for each given $\Gamma$, the associated $X$ is unique in distribution.
Consequently, $X$ is a strong Markov process (see \cite[Theorem 4.4.2]{ethi:86}).  Let ${\cal A}$ denote the set of measures $\Gamma$ for which there is some $X$ such that $(X,\Gamma)$ is a relaxed solution of the harvesting model.

A couple of observations will help the reader to understand this relaxed formulation for the model.  First, consider the solution $(X, 0$) in which the measure-valued random variable $\Gamma \equiv 0$ so it has no mass and thus no harvesting occurs.  Then Theorem~5.3.3 of \cite{ethi:86} shows the existence of a Brownian motion $W$ adapted to a possibly enlarged filtration $\{\tilde{\cal F}_t\}$ such that the process $X$ satisfies (\ref{harvested-pop}) with $Z\equiv 0$.  Next, let $Z$ denote a ``strict'' harvesting policy; that is, $Z$ is a nonnegative, increasing process that is c\`adl\`ag and adapted to $\{{\cal F}_t\}$.
Define the random measure $\Gamma$ for Borel measurable $G \subset \mathcal R$  and $t\ge 0$ by
\beq \label{gamma-def}
\Gamma(G  \times [0,t]) = \int_0^t I_{G }(X(s-), \Delta Z(s))\, dZ(s).
\eeq
It then follows that $(X,\Gamma)$ will be a relaxed solution of the harvesting model whenever $(X,Z)$ satisfies (\ref{harvested-pop}).

We turn now to the extension of the reward criterion (\ref{income-J-defn}) to the relaxed framework.  Specifically, $f: \rr_+ \mapsto \rr_+$ represents the instantaneous marginal yield accrued from harvesting. Assume $f$ is continuous and non-increasing with respect to $x$. Thus $f(x) \ge f(y)$  whenever $x\le y$; this assumption indicates that the price when the species is plentiful is smaller than when it is rare.  Moreover, we assume $0< f (0) < \infty$.  Let $(X,\Gamma)$ be a solution to the harvesting model \eqref{mgp}.  Let $S =(0,\infty)$ be the {\em survival set} of the species and denote the {\em extinction time} by $\tau = \inf\{t \geq 0: X(t) \notin S\}$.  Then the expected total discounted value from harvesting is
\begin{equation}\label{relaxed-obj-fn}
J(x_0, \Gamma) :=\ex \left[\int_{\mathcal R \times [0,{\tau}]} e^{-r s} f(x)\, \Gamma(dx  \times dz \times ds)\right].
\end{equation}
The goal is to maximize the expected total discounted value from harvesting over relaxed solutions $(X,\Gamma)$ of the harvesting model and to find an optimal harvesting strategy $\Gamma^*$.  Thus, we seek
\begin{equation}\label{value}
V(x_0) = J(x_0,\Gamma^*) := \sup_{\Gamma \in {\cal A}} J(x_0,\Gamma).
\end{equation}
We emphasize that the initial position $x_0$ is merely a parameter in the problem and that $V$ is not to be viewed as a function with any particular properties but merely is the value of the harvesting problem when the initial population size is $x_0$.  We do, however, obtain the value in functional form for $x_0$ in two regions.

\section{Linear Programming Formulation and Main Result}

Throughout this paper, we assume the equation $(\op-r) u(x)=0$ has two fundamental  solutions $\tilde\psi$ and $\tilde \phi$, where $\tilde\psi   $ is strictly increasing and $\tilde\phi$ is strictly decreasing.
As   in \cite{Alvarez}, we   put
\bed \psi(x) =\begin{cases}\tilde \psi(x), & \text{if }0 \text{ is natural or exit,}\\
\tilde\psi(x)- \frac{\tilde \psi(0)}{\tilde\phi(0)}\tilde\phi(x), & \text{if }0 \text{ is regular}. \end{cases}\eed
Note that $\psi$ solves $(\op-r) u(x)=0$, is strictly increasing, and satisfies $\psi(0)=0$.

The main result of this paper is summarized in the following theorem.

 \begin{thm}\label{thm-main} Assume that there exists  some $\tilde b \ge 0$ such that
\begin{enumerate}\item[{\em (i)}] \beq\label{key-assumption} \frac{f(x)}{\psi'(x)} \le    \frac{f(\tilde b)}{\psi'(\tilde b)},\ \ \forall x\ge 0,\eeq
 \item[{\em (ii)}]
  the function $\frac{f}{\psi'}$ is nonincreasing on $[\tilde b, \infty)$, and
 \item[{\em (iii)}] the function $f$ is continuously differentiable on $(\tilde b,\infty)$.
 \end{enumerate}
 Put $b^*=\inf\{\tilde b\ge 0: \tilde b \text{ satisfies  {\em (i)}--{\em (iii)}}\}$.
 Then the value is given by
 \beq\label{value-fn}
 V(\xz)=\begin{cases}  \displaystyle \frac{f(b^*)\psi(\xz)}{\psi'(b^*)},\ \ & \text{ if }0< \xz \le b^*, \\[2ex] \displaystyle
  \int_{b^*}^{\xz} f(y)dy + \frac{f(b^*)\psi(b^*)}{\psi'(b^*)},\ \ & \text{ if }\xz > b^*,
\end{cases}
 \eeq
 and an optimal relaxed harvesting policy is given by
 \beq\label{optimal-relaxed-policy} \Gamma^*(dx\times dz\times dt)=I_{(b^*,\infty)}(x_0) \lambda_{[b^*,\,\xz]}(dx) \delta_{\set{0}} (dz) \delta_{\set{0}}(dt) +  \Gamma_{b^*}(dx\times dz\times dt),  \eeq
where $\lambda_{[b^*,\, \xz]}(\cdot)$ denotes Lebesgue measure on $[b^*,\xz]$ and $\Gamma_{b^*}$ is defined in \propref{prop-case2}. \end{thm}

Theorem \ref{thm-main} is obtained in \cite{Alvarez} using the dynamic programming approach: the value function $V$ is obtained by explicitly solving a quasi-variational inequality of Hamilton-Jacobi-Bellman type by first using a heuristic argument to obtain $V$ and then verifying the validity of the argument.  In this paper, we use a totally different approach by imbedding the problem in a linear program over a space of measures to establish \thmref{thm-main}.  In this approach, there is no need to establish the regularity of the value function, and therefore no heuristic arguments or HJB equation are needed.  More specifically, we will first derive upper bounds for the value (depending on $x_0$), and then find a harvesting policy which achieves the appropriate upper bound.

The measures involved in the infinite-dimensional linear program are expected, discounted occupation measures corresponding to relaxed solutions $(X,\Gamma)$ of the harvesting model.  Indeed, for any Borel measurable $G_1  \subset S$ and $G \subset \mathcal R$, we define
\beq\label{measures-defn} \begin{aligned}
& \mu_\tau(G_1)= \ex\left[e^{-r\tau} I_{G_1}(X(\tau))I_{\set{\tau<\infty}}\right], \\
&\mu_0(G_1)=\ex\left[ \int_0^\tau e^{-rs} I_{G_1}(X(s))ds\right], \\
& \mu_1(G )= \ex \left[\int_{\mathcal R \times [0,\tau]} e^{-rs} I_{G}(x,z) \Gamma(dx \times  dz \times ds)\right].
\end{aligned}\eeq
Using these measures, the singular control problem of maximizing \eqref{relaxed-obj-fn} over relaxed solutions of the harvesting problem (\ref{mgp}) can be written in the form
\beq\label{lp-problem}
  \begin{cases} \text{Maximize } \disp\int f d\mu_1, \\
            \text{subject to } \disp\int g d\mu_\tau - \disp\int(\op-r)g d\mu_0 - \disp\int Bg d\mu_1 = g(x_0), \quad \forall g\in {\cal D},\\
\quad  \qquad  \qquad \mu_\tau, \mu_0, \text{ and }\mu_1
\text{ are finite measures with } \mu_\tau(S)\le 1 \text{ and }\mu_0(S)\le \frac{1}{r}.\end{cases}
\eeq

Since each relaxed solution $(X,\Gamma)$ defines measures $\mu_\tau$, $\mu_0$ and $\mu_1$ by (\ref{measures-defn}), the harvesting problem is embedded in \eqref{lp-problem}.  There might be feasible measures which do not arise in this manner.  Consequently, letting $V_{lp}(x_0)$ denotes that value of the LP problem \eqref{lp-problem} with initial condition $X(0-)=x_0 > 0$, we have
 \beq\label{value-ineq1}
 V(x_0) \le V_{lp}(x_0).
 \eeq

\section{The Proof of \thmref{thm-main}}
This section is devoted to the proof of \thmref{thm-main}.  We consider two different cases: when $0< x_0 \le b^*$ and when $ x_0 > b^*$, where $b^*$ is the threshold level given in the statement of the theorem. 

\subsection{Case 1: $0< x_0 \le b^*$}
Our goal is to find the value $V(x_0)$ defined in \eqref{value} and a relaxed optimal harvesting policy directly. The proof follows along the lines of the arguments used in \cite{helm:11}.  In fact, an optimal strict harvesting policy $Z$ is obtained so the relaxed formulation is not necessary in this case.  The general argument involves finding an upper bound for $V_{lp}(x_0)$ by reducing the number of constraints in the linear program \eqref{lp-problem} and then identifying a solution $(X^*,Z^*)$ which achieves the bound.  The relaxed harvesting policy $\Gamma^*$ is obtained from $Z^*$ by \eqref{gamma-def}.

We will need the Skorohod lemma (see \cite{Lions-S-84}) so we give its statement for completeness.
\begin{lem}\label{lem-Skorohod}
Given any initial state $\xz$ and any boundary $c$, there exists a unique $\{{\F}_t\}$-adapted c\`{a}dl\`{a}g pair  $(X, L_c)$ such that $L_c$ is nonnegative and nondecreasing and
\begin{align}
\label{Skorohod-eq1}&  X(t)= \xz + \int_0^t b(X(s))ds + \int_0^t\sigma(X(s))dW(s) - L_c(t), \\
\label{Skorohod-eq2}& X(t) \in (-\infty,c], \ \text{ for almost all }\  t\ge 0, \\
\label{Skorohod-eq3}& \int_0^\infty I_{\set{X(s)< c}}dL_c(s) =0.
\end{align}
Moreover, $L_c$ is continuous if $\xz  \le  c$.
\end{lem}
The solution $X$ to the above equations is a reflected diffusion at the boundary $c$, and the process $L_c$ is the local time process of $X$ at $c$. Moveover, the property \eqref{Skorohod-eq3} shows that the process $L_c$ increases only when $X$ reaches the boundary $c$.

\begin{prop}\label{prop-case1} Let $0< \xz \le b^*$.  Then
\bed V(\xz)= \frac{f(b^*)\psi(\xz)}{\psi'(b^*) } , \eed
and an optimal harvesting strategy is given by the local time process $L_{b^*}$ of $X^*$ at $b^*$.
\end{prop}

\begin{proof}
Though $\psi$ does not have compact support, an argument similar to the one in
 \cite{helm:11}
 shows we may use the function $\psi$ in the constraints of \eqref{lp-problem}.  This results in an auxiliary linear program
\beq\label{auxiliary-lp}
 \begin{cases} \text{Maximize } \disp\int f d\mu_1, \\
            \text{subject to }
  - \disp\int B\psi d\mu_1 = \psi(x_0),  \\ \hfill
        \phantom{subject to }    \mu_1
\text{ is a finite measure.} \end{cases}
\eeq
In obtaining the auxiliary linear program we have used the properties that $\psi(0)=0$ and $(A-r)\psi(x)=0$ to eliminate the measures $\mu_\tau$ and $\mu_0$ from the program.  Denote  the solution to  \eqref{auxiliary-lp} by $V_{aux} (x_0)$. Then  since \eqref{auxiliary-lp} has fewer constraints than \eqref{lp-problem}, the set of feasible measures $\mu_1$ for \eqref{auxiliary-lp} may contain more $\mu_1$ measures than those arising from the feasible solutions to \eqref{lp-problem} and hence
\beq\label{value-ineq2} V(\xz)\le V_{lp}(\xz)\le V_{aux}(\xz).\eeq

Using the definition of $B$ in \eqref{genB},  the constraint in \eqref{auxiliary-lp} can be written as
\bed \psi(x_0) = - \int B\psi d\mu_1 =  \int_{\cal R} \(\psi'(x)I_{\set{0}}(z) + \frac{\psi(x)-\psi(x-z)}{z}I_{(0,x]}(z) \)\; \mu_1(dx\times dz).\eed
Recall that $\psi$ is strictly increasing and $\psi(0)=0$. Therefore $\psi(\xz)>0$ and hence it follows that
\bed 1=  \int_{\cal R} \frac{\psi'(x)I_{\set{0}}(z) + \frac{\psi(x)-\psi(x-z)}{z}I_{(0,x]}(z) }{\psi(\xz)} \mu_1(dx\times dz). \eed
Thus the integrand is a probability density relative to any feasible measure $\mu_1$ and defines a corresponding probability measure $\tilde \mu_1$ on ${\cal R}$.
Now the objective function \eqref{relaxed-obj-fn} can be rewritten as
\beq \label{obj-fn-2} \int fd\mu_1 =  \int \frac{f(x) \psi(\xz)}{\psi'(x)I_{\set{0}}(z) + \frac{\psi(x)-\psi(x-z)}{z}I_{(0,x]}(z)}\;  \tilde\mu_1(dx\times dz).
\eeq

We claim that
\beq\label{f/psi-ineq}
\frac{f(x)  }{\psi'(x)I_{\set{0}}(z) + \frac{\psi(x)-\psi(x-z)}{z}I_{(0,x]}(z)} \le \frac{f(b^*)}{\psi'(b^*)}.
\eeq
In fact for $z=0$, \eqref{key-assumption} and the definition of $b^*$ in the statement of Theorem~\ref{thm-main} implies
\bed \frac{f(x)  }{\psi'(x)I_{\set{0}}(z) + \frac{\psi(x)-\psi(x-z)}{z}I_{(0,x]}(z)}  = \frac{f(x)}{\psi'(x)} \le \frac{f(b^*)}{\psi'(b^*)}. \eed
On the other hand, for $z\not= 0$, then the assumption that $f$ is nonincreasing along with \eqref{key-assumption} implies that for some $\theta \in [0,1]$
\begin{eqnarray*}
\frac{f(x)  }{\psi'(x)I_{\set{0}}(z) + \frac{\psi(x)-\psi(x-z)}{z}I_{(0,x]}(z)} &=& \frac{f(x)z}{\psi(x)-\psi(x-z)} = \frac{f(x)z}{\psi'(x-\theta z) z} \\
&=& \frac{f(x-\theta z)}{\psi'(x-\theta z)} \frac{f(x)}{f(x-\theta z)} \le \frac{f(x-\theta z)}{\psi'(x-\theta z)} \le \frac{f(b^*)}{\psi'(b^*)}.
\end{eqnarray*}

Now it follows from \eqref{obj-fn-2} and the bound in \eqref{f/psi-ineq} that for any feasible measure $\mu_1$ of \eqref{auxiliary-lp}
 \bed \int fd\mu_1 \le  \int_{{\cal R}}\frac{f(b^*)}{\psi'(b^*)} \psi(\xz) \tilde\mu_1(dx\times dz) \le \frac{f(b^*)}{\psi'(b^*)} \psi(\xz),\eed
and hence \beq\label{case1-ineq} V_{aux}(\xz) \le \frac{f(b^*)}{\psi'(b^*)} \psi(\xz). \eeq

Next we show that there is an admissible (strict) harvesting strategy $Z^*$ and therefore a relaxed harvesting strategy $\Gamma^* \in \cal A$ such that
\beq\label{case1-claim} J(x_0,Z^*) = J(\xz,\Gamma^*)= \frac{f(b^*)}{\psi'(b^*)} \psi(\xz).\eeq

Recall, we are analyzing the case in which $\xz \le b^*$.  Let $(X^*, L_{b^*})$ be the solution to the Skorohod problem \eqref{Skorohod-eq1}--\eqref{Skorohod-eq3} of \lemref{lem-Skorohod} with $c = b^*$.  Note that $L_{b^*}$ is continuous and hence $X^*$ is also continuous.  Next for any $t>0$, by virtue of It\^o's formula and \eqref{Skorohod-eq3}, we have
 \beq\label{eq-ito} \begin{aligned} \ex_{\xz} & [e^{-r (\tau \wedge t)} \psi(X^*(\tau \wedge t))] -\psi(\xz)\\
 &  = \ex_{\xz} \left[\int_0^{\tau \wedge t} e^{-rs} (\op-r) \psi(X^*(s)) ds
 -  \int_0^{\tau \wedge t} e^{-rs} \psi'( X^*(s)) dL_{b^*}(s)\right] \\
 & =- \psi'(b^*) \ex_{\xz}\left[ \int_0^{\tau \wedge t} e^{-rs}   dL_{b^*}(s)\right].  \end{aligned}\eeq
Due to the process $X^*$ being bounded (from \eqref{Skorohod-eq2}), $\psi(X^*(t))$ is also bounded for all $t \geq 0$.  This observation along with the fact that $\psi(0)=0$ then implies
\bea \ad \lim_{t\to \infty} \ex_{\xz}  [e^{-r (\tau \wedge t)} \psi(X^*(\tau \wedge t))] \\ \aad \ = \lim_{t\to\infty}  \ex_{\xz}  \left[e^{-r  t} \psi(X^*( t)) I_{\set{\tau=\infty}} + e^{-r (\tau \wedge t)} \psi(X^*(\tau \wedge t)) I_{\set{\tau < \infty}}\right] =0.  \eea
Hence by letting $t \to \infty$  in \eqref{eq-ito}, it follows that
\bed \ex_{\xz} \left[\int_0^{\tau} e^{-rs}   dL_{b^*}  (s)\right]= \frac{\psi(\xz)}{\psi'(b^*) }, \eed
which in turn implies that
$$J(\xz, L_{b^*}) = \ex_{\xz} \left[\int_0^\tau  e^{-rs}  f(X^*(s))dL_{b^*}(s)\right] = f(b^*) \ex_{\xz}\left[ \int_0^\tau  e^{-rs}  dL_{b^*}(s)\right] = \frac{f(b^*)\psi(\xz)}{\psi'(b^*) }.$$
Therefore \eqref{case1-claim} follows with $Z^*= L_{b^*}$. Defining $\Gamma^*$ by \eqref{gamma-def}, the pair $(X^*,\Gamma^*)$ is a relaxed solution of the harvesting model which achieves the bound.
\end{proof}

Since $\Delta L_{b^*}(s)=0$ for every $s\ge 0$, an optimal strategy is to harvest just enough of the population (using the local time of $X^*$ at $b^*$) so that the population size ``reflects'' at $b^*$.

\subsection{Case 2: $  x_0 > b^*$}

This case is the more interesting of the two cases and requires a new argument and also a different type of harvesting policy than what appears in the literature.  It is for this case that the relaxed formulation of the problem is needed in order to obtain an optimal control.

When dealing with singular control problems, one usually takes the so-called reflection strategy, namely,
\beq\label{reflection-strategy}Z(t)= (\xz-b^*)^+ + L_{b^*}(t), \eeq
where one uses the local time process $L_{b^*}$ at $b^*$ following an immediate jump from $x_0$ to $b^*$.  Such a reflection strategy is used in  \cite{choulli}, \cite{Pha09} and others.  The income corresponding to \eqref{reflection-strategy} is
$$J(\xz, Z)= f(\xz)(\xz-b^*) +  \frac{f(b^*)\psi(b^*)}{\psi'(b^*) }.$$
When $f$ is strictly decreasing, then the reflection strategy is not optimal. In fact, there is no  strict admissible optimal harvesting strategy; please see \cite{Song-S-Z} for detailed arguments as well as the explicit construction of an $\e$-optimal admissible harvesting policy for a regime-switching diffusion (the static environment model of this paper being a special case).

Our purpose is to find an optimal relaxed harvesting strategy. The previous section proves $V(\xz) \le \frac{f(b^*)\psi(\xz)}{\psi'(b^*) }$. However, the upper bound is a strict upper bound; no relaxed harvesting policy will achieve this upper bound. The following arguments determine a sharp upper bound.  We begin by establishing the following estimate.
\begin{lem}\label{lem2-case2}
Assume the conditions in \thmref{thm-main}. Denote \beq\label{g-fn-defn}
g(x):= \int_{b^*}^x f(y)dy, \ \ \text{ for } x \ge b^*. \eeq
Then
\beq\label{op-rg<?}
(\op-r) g(x) \le r \frac{f(b^*)\psi(b^*)}{\psi'(b^*)}, \ \ \text{for every }x > b^*.
\eeq
\end{lem}
\begin{proof}
Since the function $f/\psi'$ is nonincreasing on $(b^*,\infty)$, we have
\bed 0\ge \frac{d}{dx}  \(\frac{f(x)}{\psi'(x)}\) = \frac{f'(x)\psi'(x)-f(x)\psi''(x)}{(\psi'(x))^2}, \ \ x > b^* . \eed
But $\psi$ is strictly increasing and so $\psi'(x)>0$. Hence it follows that  $ f'(x)\psi'(x)-f(x)\psi''(x) \le 0$,   or equivalently $$ f'(x) \le \frac{f(x)}{\psi'(x)}\psi''(x), \text{ for }x> b^*. $$
It then follows that for each $x > b^*$
\bed \begin{aligned}
(\op-r) g(x) &
  = \frac{1}{2} \sigma^2(x) f'(x) + b(x) f(x) - r \int_{b^*}^{x}\frac{f(y)}{\psi'(y)} \psi'(y)dy \\
& \le \frac{1}{2} \sigma^2(x) f'(x) + b(x) f(x) - r \int_{b^*}^{x}\frac{f(x)}{\psi'(x)} \psi'(y)dy\\
&\le  \frac{1}{2} \sigma^2(x) \frac{f(x)}{\psi'(x)} \psi''(x) + b(x) f(x) - r\frac{f(x)}{\psi'(x)} (\psi(x)- \psi(b^*))\\
& = \frac{f(x)}{\psi'(x)} \left[\frac{1}{2} \sigma^2(x) \psi''(x)  + b(x) \psi'(x) -r \psi(x)\right] + r\frac{f(x)}{\psi'(x)}\psi(b^*)\\
&=r\frac{f(x)}{\psi'(x)}\psi(b^*) \le r\frac{f(b^*)}{\psi'(b^*)}\psi(b^*).
\end{aligned} \eed
 \end{proof}

The next result establishes a sharper upper bound on the value of the problem.  This upper bound will be seen to be the value of the harvest for a relaxed solution of the harvesting model and hence establishes the value.

\begin{prop}\label{prop4-case2}
Let $x_0 > b^*$ and assume the conditions of \thmref{thm-main}.  Then
\beq\label{payoff-upper-bd}
V(\xz) \le \int_{b^*}^{\xz} f(y)dy + \frac{f(b^*)\psi(b^*)}{\psi'(b^*)}.
\eeq
\end{prop}

\begin{proof} Let $(X,\Gamma)$ be an arbitrary solution to the harvesting model \eqref{mgp} and define
$$\tau_{b^*}= \inf\{t>0: X(t)\le b^*\}$$
and observe that $\taub\le \tau$.  The rest of the proof is divided into several steps.
\medskip

\noindent
{\em Step 1.}\/  We claim that
\beq\label{case2-claim1} J(\xz, \Gamma)
\le \ex_{\xz}\left[\int_{\R\times [0,\taub]}
 e^{-rs} f(x)\Gamma(dx\times dz\times dt)\right]+ \frac{f(b^*)\psi(b^*)}{\psi'(b^*)} \ex_{\xz} \left[e^{-r\taub} \right].\eeq
 To establish \eqref{case2-claim1}, we write
 \beq\label{eq-J-I+II} \begin{aligned}
 J(\xz, \Gamma) & = \ex_{\xz} \left[\int_{\mathcal R \times [0,\tau]} e^{-r s} f(x) \Gamma(dx\times dz\times ds)\right]\\
  & = \ex_{\xz} \left[\int_{\mathcal R \times [0,\tau]} e^{-r s} f(x) \Gamma(dx\times dz\times ds)(I_{\set{\taub=\infty}}+I_{\set{\taub<\infty}})\right]. 
 \end{aligned}\eeq
Clearly on the set $\{\taub = \infty\}$ we also have $\tau = \infty = \taub$ so the first term can be rewritten as
\beq\label{eq-J-I} 
\ex_{\xz} \left[\int_{\mathcal R \times [0,\taub]} e^{-r s} f(x) \Gamma(dx\times dz\times ds) I_{\set{\taub=\infty}}\right]. 
\eeq
For the second term,  it follows from the strong Markov property and \eqref{case1-ineq} that
\begin{eqnarray*}
\lefteqn{\ex_{\xz}\left[\int_{\R \times [0,\tau] } e^{-rs} f(x) \Gamma(dx\times dz \times ds)I_{\set{\taub<\infty}}\right]} \\
&=& \ex_{\xz}\bigg[ I_{\set{\taub<\infty}}\bigg( \int_{\R\times [0,\taub)} e^{-rs} f(x)\Gamma(dx\times dz \times ds) \\
& & \qquad \qquad\qquad\quad +  \ex_{\xz}\bigg[\ex_{\xz} \left[\int_{\R\times [\taub,\tau]} e^{-rs} f(x)\Gamma(dx\times dz\times ds)\bigg|\F_{\taub}
\right] \bigg]\bigg)\bigg] \\
&\le& \ex_{\xz}\bigg[ I_{\set{\taub<\infty}}\bigg( \int_{\R\times [0,\taub]} e^{-rs} f(x)\Gamma(dx\times dz \times ds) \\
& & \qquad \qquad\qquad\quad+  e^{-r\taub}\ex_{ X(\taub)} \left[\int_{\R\times [0, \tau]} e^{-rs} f(x)\tilde\Gamma(dx\times dz\times ds) \right]\bigg)\bigg] \\
&\le& \ex_{\xz}\left[I_{\set{\taub<\infty}} \(\int_{\R\times [0,\taub]}
 e^{-rs} f(x)\Gamma(dx\times dz\times ds)+ e^{-r\taub} \frac{f(b^*)}{\psi'(b^*)}\psi(X(\taub))\)\right],
\end{eqnarray*}
where $\tilde \Gamma(G\times[0,t])= \Gamma(G\times[\taub, \taub+t])$ for $G\subset \R$.  But on the set $\set{\taub<\infty}$, $X(\taub) \le b^*$. Note also that $\psi $ is strictly increasing. Thus we have
\begin{equation} \label{eq-J-II}  \begin{array}{l} \displaystyle
\ex_{\xz}\left[\int_{\R \times [0,\tau] } e^{-rs} f(x) \Gamma(dx\times dz \times ds)I_{\set{\taub<\infty}}\right] \rule[-15pt]{0pt}{15pt}\\ \displaystyle
\quad \le  \ex_{\xz}\left[\int_{\R\times [0,\taub]}
 e^{-rs} f(x)\Gamma(dx\times dz\times dt)I_{\set{\taub<\infty}}\right]+  \frac{f(b^*)\psi(b^*)}{\psi'(b^*)} \ex_{\xz} \left[e^{-r\taub} \right].
 \end{array}
\end{equation}
Finally a combination of \eqref{eq-J-I+II}--\eqref{eq-J-II} implies \eqref{case2-claim1}.
\medskip

\noindent
{\em Step 2.} Since $f$ is nonincreasing, for any $x,\delta > 0$ with $x-\delta \ge b^*$, we have
\bed f(x) \delta \le \int_{x-\delta}^x f(y) dy = g(x)-g(x-\delta). \eed
Therefore it follows that
\begin{eqnarray*}
\lefteqn{\ex_{\xz}\left[ \int_{\R\times [0, \taub]} e^{-rs}f(x) \Gamma(dx\times dz\times dt)\right]} \\
&=&  \ex_{\xz} \left[  \int_{\rr_+\times\set{0}\times [0, \taub]} e^{-rs}f(x) \Gamma(dx\times dz\times dt)\right] \\
& & \quad  + \ex_{\xz}\left[   \int_{(\R-(\rr_+\times\set{0}))\times [0, \taub]} e^{-rs}\frac{f(x)z}{z} \Gamma(dx\times dz\times dt)  \right] \\
&\le&   \ex_{\xz} \left[  \int_{\rr_+\times\set{0}\times [0, \taub]} e^{-rs}f(x) \Gamma(dx\times dx\times dt)\right] \\
& & \quad +    \ex_{\xz}\left[   \int_{(\R-(\rr_+\times\set{0}))\times [0, \taub]} e^{-rs}\frac{g(x)-g(x-z)}{z} \Gamma(dx\times dx\times dt)  \right].
\end{eqnarray*}
Recalling the definition of $B$ in \eqref{genB}, we observe that $Bg(x,0) = -g'(x) = - f(x)$ and $\frac{g(x) - g(x-z)}{z} = -Bg(x,z)$ when $z > 0$ and hence
\beq\label{case2-claim2} \begin{array}{l} \displaystyle
\ex_{\xz} \left[\int_{\R\times [0, \taub]} e^{-rs}f(x) \Gamma(dx\times dx\times dt)\right] \\ \displaystyle
\qquad \qquad \le  - \ex_{\xz} \left[\int_{\R\times [0,\taub]} e^{-rs} Bg(x,z)\Gamma(dx\times dx\times dt)\right].
\end{array}
\eeq
\smallskip

\noindent
{\em Step 3.}\/  We have
\beq\label{E-Bg<=} - \ex_{\xz}\left[ \int_{\R\times[0,\taub]} e^{-rs} Bg(x,z)\Gamma(dx\times dz\times dt)\right] \le g(\xz) + \frac{f(b^*)\psi(b^*)}{\psi'(b^*)} (1-\ex_{\xz}[e^{-r \taub}] ). \eeq
 In fact, for any $t>0$,
It\^o's formula  implies  that
\bed\begin{aligned} \ex_{\xz}& [e^{-r (\taub\wedge t)}g(  X(\taub\wedge t))] - g(\xz) \\ &= \ex_{\xz} \left[  \int_{0}^{\taub\wedge t} e^{-rs} (\op -r )g(  X(s))ds
+ \int_{\R\times [0, \taub\wedge t]} e^{-rs}  Bg(x,z)\Gamma(dx\times dx\times dt)\right].\end{aligned}\eed
Isolating the term involving $Bg$ and using the bound \eqref{op-rg<?}, we have
\bed\begin{aligned} - & \ex_{\xz}  \left[  \int_{\R\times [0, \taub\wedge t]} e^{-rs}  Bg(x,z)\Gamma(dx\times dx\times dt)\right]
\\ & \le g(\xz)
 -\ex_{\xz}  [e^{-r (\taub\wedge t)}g(  X(\taub\wedge t))]
  + \ex_{\xz} \left[  \int_0^{\taub\wedge t} e^{-rs} r \frac{f(b^*)\psi(b^*)}{\psi'(b^*)} ds\right]\\
& \le g(\xz) + \frac{f(b^*)\psi(b^*)}{\psi'(b^*)} (1- \ex_{\xz}[e^{-r(\taub\wedge t)}]). \end{aligned}\eed
Now \eqref{E-Bg<=} follows by letting $t \to \infty$ in the above inequality.
\medskip

\noindent
{\em Step 4.}\/  Combining  \eqref{case2-claim1},  \eqref{case2-claim2},  and \eqref{E-Bg<=} yields
\bed \begin{aligned}J(\xz, \Gamma) & \le g(\xz) + \frac{f(b^*)\psi(b^*)}{\psi'(b^*)} (1-\ex_{\xz}[e^{-r \taub}] ) + \frac{f(b^*)\psi(b^*)}{\psi'(b^*)} \ex_{\xz} [e^{-r\taub} ] \\
 & = g(\xz) + \frac{f(b^*)\psi(b^*)}{\psi'(b^*)} .
\end{aligned}\eed
The bound in \eqref{payoff-upper-bd} is therefore established by taking supremum over $\Gamma \in {\cal A}$.
\end{proof}

We have derived an upper bound for the value $V(\xz)$ in \propref{prop4-case2}.  The next natural question is: ``Can we find an admissible optimal harvesting policy which  achieves the upper bound specified in the right-hand side of \eqref{payoff-upper-bd}?''  The following proposition answers this question in the affirmative by explicitly constructing an optimal relaxed harvesting policy.

\begin{prop}\label{prop-case2}
Let $\lambda_{[b^*,  \xz]}(\cdot)$ denote Lebesgue measure on $[b^*, x]$.  Also let $L_{b^*}$ denote the local time process of \propref{prop-case1} with $x_0$ taken to be $b^*$.  Define $\Gamma_{b^*}$ to be the random measure defined by \eqref{gamma-def} using $Z=L_{b^*}$.  Finally, define the relaxed harvesting strategy by
$$ \Gamma^*(dx\times dz\times dt)= \lambda_{[b^*,  \xz]}(dx) \delta_{\set{0}} (dz) \delta_{\set{0}}(dt) + \Gamma_{b^*}(dx\times dz\times dt).$$
Then
\beq \label{equality}
V(\xz)=J(\xz, \Gamma^*) =\int_{b^*}^{\xz} f(y)dy + \frac{f(b^*)\psi(b^*)}{\psi'(b^*)}.
\eeq
\end{prop}

\begin{proof}
We observe that the measure $\mu_1$ obtained from $\Gamma^*$ by \eqref{measures-defn} is
$$\mu_1^*(dx\times dz)= \left[\lambda_{[b^*,\xz]}(dx) + \frac{\psi(b^*)}{\psi'(b^*)}\delta_{\set{b^*}}(dx)\right]\times \delta_{\set{0}}(dz)$$
and is feasible for \eqref{auxiliary-lp}.  The measure $\lambda_{[b^*,  \xz]}(\cdot) \delta_{\set{0}} (\cdot) \delta_{\set{0}}(\cdot)$ instantaneously resets the problem at time $0$ so that the initial position of the population becomes $b^*$.  Take $(X^*,L_{b^*})$ to be the solution of the Skorhod problem from \lemref{lem-Skorohod} with $x_0=b^*$.  It is then easy to verify that $(X^*,\Gamma^*)$ is a relaxed solution to the harvesting model whose value equals the right-hand side of \eqref{equality}.
\end{proof}

We observe that the manner in which this optimal harvesting policy differs from the typical ``reflection'' strategy occurs at the initial time.  Whereas the reflection strategy has the process $X$ instantaneously jump from $x_0$ to $b^*$, the optimal relaxed harvesting policy obtains this relocation in an instantaneous {\em but continuous}\/ manner.
\medskip

Finally we note that the combination of Propositions \ref{prop-case1} and \ref{prop-case2} establishes \thmref{thm-main}.  Moreover, the optimal relaxed harvesting policy can be written as
 $$ \Gamma^*(dx\times dz\times dt)= I_{(b^*,\infty)}(x_0) \lambda_{[b^*,  \xz]}(dx) \delta_{\set{0}} (dz) \delta_{\set{0}}(dt) + \Gamma_{b^*}(dx\times dz\times dt),$$
which unifies the two cases.


\def\cprime{$'$}

\end{document}